\documentclass[11pt]{article}

\usepackage{amsmath, amssymb, amsthm, amsfonts, dsfont, bbm,mathtools}
\usepackage{enumerate}
\usepackage{a4wide}
\usepackage[table]{xcolor}
\usepackage{nicefrac}
\usepackage{csquotes}

 \usepackage{natbib}
\def\phi{\varphi}

\mathchardef\sim="2218

\newcommand{\atm}{\ensuremath{\mathrm{AST}}}
\newcommand{\zf}{\ensuremath{\mathrm{ZF}}}
\newcommand{\zfc}{\ensuremath{\mathrm{ZFC}}}
\newcommand{\zffin}{\ensuremath{\mathrm{ZF}_{\mathrm{fin}}}} 
\newcommand{\gb}{\ensuremath{\mathrm{GB}}}

\newcommand{\tss}{\ensuremath{\mathrm{TSS}}}

\newcommand{\pa}{\ensuremath{\mathrm{PA}}}

\newcommand{\FN}{\ensuremath{\mathrm{FN}}} 
\newcommand{\N}{\ensuremath{\mathrm{N}}}

\title{Vopěnka's Alternative Set Theory as a framework for \\feasible numbers}
\author{Zuzana Haniková\\
Institute of Computer Science of the Czech Academy of Sciences
}
   
\begin{document}

\maketitle

\begin{abstract}
Vopěnka's Alternative Set Theory has been considered as a framework for modelling vague notions. 
This paper takes feasibility, pertaining to numbers as per some of Yessenin-Volpin's work, and tries to assess how this notion could be modelled in the Alternative Set Theory. The route explored in detail consists in an attempt to model Dummett's weakly infinite, weakly finite totalities, the coherence of which Dummett takes to be decisive of the coherence of Yessenin-Volpin's foundational conception in general. The outcome of the particular analysis turns out to be negative: the Alternative Set Theory could be taken as a model of some, but not quite all features of Yessenin-Volpin's feasible numbers. Given that only one particular approach is explored here, the outcome has limited bearing on the more general questions of coherence of ultrafinitist position, or indeed the possibilities of modelling feasible numbers.
\end{abstract}

\section{Introduction}

I take \emph{ultrafinitism} and \emph{strict finitism} to be  synonymous, or at least synonymous enough to try to find some common ground between the various discussions that go under their various respective names. Along with \emph{ultraintuitionism} or \emph{anthropologism},\footnote{This term, coming from  \cite{WangHao:80years}, is explicitly stated to be a \textquote{more colorful} one for strict finitism earlier discussed by \cite{Kreisel:reviewWittgRFM}.} these names and the preferences for one or another may be due to the various language backgrounds of the proponents.

Ultrafinitism has been  deliberated as a foundational position: notably,
Yessenin-Volpin highlighted the revisionist facet of ultrafinitism in his proposed programme 
\cite{Yessenin-Volpin:ProgrammeUltra, Yessenin-Volpin:UltraIntCriticism1970}. 
Long before that, Bernays in the  lecture published as \cite{Bernays:OnPlatonism} 
pointed out that intuitionism as a foundational position, 
which stands perhaps as a typical example of a revisionist programme (even if not all proponents of 
constructive or intuitionist mathematics take the revisionist view), might  be found unstable if its tenets are 
taken literally; or if one swaps evidence obtainable in principle for evidence obtained in practice. 
These remarks of Bernays are not to be taken as an attempt to buttress the case for ultrafinitism or to undermine
intuitionism; one  takes away rather the impression that Bernays took a pluralistic position. 
While the lecture does not aim to strongly defend any singular foundational view,  
it rebukes the revisionist attitude of the intuitionists and  tries to show the merit of having both the realist and the intuitionist position on the table.  
As pointed out in  \cite{Hallett:CantorianSetTheory},
long before Bernays in turn, Cantor would argue for his (realist) position about sets by pointing out the 
practical limits of our capacity of generating and evidencing the natural numbers as large but finite sets.

These remarks aim  to indicate why one might take an interest in the matter: Bernays pointed out the constructive direction
of departure from the mainstream mathematics (he refers to the latter as \textquote{platonist})\footnote{\cite{WangHao:80years} remarks that the position exemplifies
a \textquote{mathematics of being}.}  and raised the question how far
a mathematician might be comfortable with departing in this particular direction.

This paper takes a notion that epitomizes  Yessenin-Volpin's work, namely \textquote{feasible number}
as used in \cite{Yessenin-Volpin:ProgrammeUltra, Yessenin-Volpin:UltraIntCriticism1970},
and tries to model it using the notions native to the Alternative Set Theory, as developed by
Vopěnka and a group of researchers he worked with and presented first in a comprehensive form in the monograph \cite{Vopenka:Teubner}.
The paper is written  from a perspective of familiarity with Vopěnka's work.
A question that is of interest to me is in which sense and to which extent Vopěnka's work might be 
understood as a contribution to ultrafinitistic mathematics. I am not planning to answer this broader question here in full; 
the aim is to analyze the issue of feasible numbers  on the backdrop of the Alternative Set Theory.

Linking  Yessenin-Volpin and Vopěnka is by no means a novelty: Vopěnka's book \cite{Vopenka:Teubner}  refers to one of Yessenin-Volpin's works in Russian.  
Another influence on Vopěnka that needs to be taken into account, although
it can at present be accounted for rather less comprehensively than one might wish,
is that of Rieger, whose seminars Vopěnka frequented as a young researcher.\footnote{When Rieger passed away in 1963, Vopěnka started his own seminar in axiomatic set theory, cf.~\cite{Vopenka:PragueSetTheorySeminar}.} 
Rieger devoted his research   to mathematical logic and foundational topics following his stay in Warsaw with the Mostowski group,\footnote{Cf.~\cite{Pecinova:Rieger}.}
and participated\footnote{Cf.~\cite{Rieger:SurProblemeNaturels}.} in the \emph{Infinitistic methods} congress in Warsaw in 1959, where 
Yessenin-Volpin was also a participant, see \cite{Yessenin-Volpin:ProgrammeUltra}.\footnote{Yessenin-Volpin in 
turn refers to Rieger in \cite{Yessenin-Volpin:UltraIntCriticism1970}, along with other names who are all reported to be 
sceptical about a categorical view of natural numbers.}

\section{Yessenin-Volpin's views on feasible numbers}

It will be useful to refer to several passages from the opening of \cite{Yessenin-Volpin:UltraIntCriticism1970} here.  

\begin{displayquote}
Let us consider the series $F$ of feasible numbers, i.e., of those up to which
it is possible to count. The number $0$ is feasible and if $n$ is feasible then
$n < 10^{12}$ and so $n^\prime$ also is feasible. 
And each feasible number can be obtained from $0$ by adding ${}^\prime$; 
so $F$ forms a natural number series. But $10^{12}$ does not belong to $F$. (p.~5)
\end{displayquote}

Two remarks on this. Yessenin-Volpin says that feasibility is an ``empirical notion''.
So in particular, it is not an arithmetical predicate, so induction would not  
necessarily apply to it.
This key fact is reflected in the seminal mathematical study of feasibility provided in
\cite{Parikh:Feasibility}. 

The other remarkable thing is the notion of a natural number series. 
The indefinite article is indicative of the approach Yessenin-Volpin takes, admitting for
consideration not just one natural number series, but rather several such.
He states his position quite explicitly by listing the assumption \textquote{T1. The uniqueness (up to isomorphism) of the natural number series} among a list of assumptions he
plans to criticize \textquote{simply on the ground they are
assumptions} (of the established, classical, \textquote{traditional} mathematics).
So while  $F$ denotes the series of feasible numbers, 
$N$ denotes a natural number series with $10^{12}$ in it.

\begin{displayquote}
It may seem quite natural to say that some numbers of $N$ are feasible and
that the number $10^{12}$ is not. I said this above myself. But this way of considering the numbers of $F$ leads to the situation that the finite set of numbers $\{0,1, \dots, 10^{12})$ of $N$ contains an infinite part $F$. Perhaps it is not absurd but it leads to many difficulties. If somebody looks over this 
set of $N$-numbers he exhausts this infinite part, but I prefer to say an infinite process
cannot be exhausted. This difficulty disappears if we distinguish between the $F$-numbers and the $N$-numbers equivalent to them. But sometimes it will be more convenient to identify these numbers belonging to different series. It is a very common thing; sometimes we identify objects that we distinguish
in other cases. The rules of doing so must be established in a rigorous way. (p.~7)
\end{displayquote}

Yessenin-Volpin's works do not adhere to the usual scruple of axiomatic mathematics---indeed perhaps one of its main features---of a very restricted vocabulary. The exposition relies on a voluminous array of notions. 
Some of these notions may be of metaphorical nature. 
However there are important notions  central to his theory
that are not metaphorical (contrary to classical mathematics, where they are almost invariably metaphorical):
the terminology that reflects the passage of time. 
Needless to say, this is so by design, as 
Yessenin-Volpin acknowledges intuitionism as the direction of his departure from classical mathematics:
\textquote{I accept the traditional intuitionistic criticism of Brouwer and go further.} 
Hence the name that he selected for his revision of the foundations, \emph{ultra-intuitionism}.

The central notion of feasibility arrives on the scene while retaining 
all its commonplace sense of action. 
So Yessenin-Volpin's analysis will pertain to what Wang 
called the \textquote{mathematics of doing}.
One can go on and provide the numerous references to \emph{constructions} of various kinds, 
to the concept of \emph{exhaustion}, to 
a purported theory of \emph{collations} (identifications and distinctions),
to Zenonian situations, 
the talk of tactics, actions and their possibility, and perhaps also the liberalism/despotism polarity,
as examples of this phenomenon. Let me offer another quotation.

\begin{displayquote}
In virtue of the traditional intuitionistic criticism we have to consider
the natural number series not as accomplished totalities but rather as
infinite processes.
But one cannot speak exactly about any process without
using tenses in order to discern the past events from the future ones. For
a natural number series, at each stage, its infiniteness belongs to the domain
of the future and so one has to expose the rules of transforming the future
tense into the present or the past in reasonings. (p.~12)
\end{displayquote}
It is in the light of this quotation that I would like briefly to discuss Yessenin-Volpin's
notion of infinity, so frequently employed in the text. 

Because of the indicated (ultra-)intuitionistic setting, it cannot be taken  for granted that \textquote{finite} and \textquote{infinite}
are complementary notions in the usual (i.e., classical) sense, i.e.: wherever one 
of these terms can be used meaningfully, there the other one can be used too; and 
each item in the domain of their application is finite, or is infinite, but never both at the same time.
On inspection, \emph{infinite} is a  metamathematical notion for Yessenin-Volpin;
we encounter infinity as an attribute of a series (in saying that Yessenin-Volpin rejects the categorical view of natural numbers, the notion of a \textquote{series} is interpreted as  a \textquote{structure}) or a process, 
namely one that \textquote{cannot be exhausted}. 

It is not possible to  count up to $10^{12}$.
But there are other examples in the texts that suggest 
that we should understand a series to be a process of occurrence or perhaps creation
of their members (this proceeds in \textquote{steps}), so in particular, the process can be
 in progress as one speaks. 
Then the occurrence of members are \textquote{events}.
It is possible for  a series to be rather short: \textquote{let $N_1$ be shorter than $b_{20}$, 
which means that there never shall be such entity as $a_{20}$ in $N_1$.}
It matters that the numbers in a series differ from one another:
\textquote{In all cases the events, as they appear, are required
to be distinguished from all preceding ones; and the fulfilment of this
requirement is to be given by a single action with a parameter m for the
preceding numbers.} (p.~26)
So one could not obtain a series by always supplying the same member.

In the second quotation, its author admits that it \textquote{leads to many difficulties} 
to accommodate an infinite series, such as $F$, in the finite set  of numbers $\{0,1, \dots, 10^{12} \}$. 
His point is to highlight that an identification of the $F$-numbers with the $N$-numbers, so to speak,
is an option that one might  take (or not) and if one does, then one should be explicit about it and 
possibly supply some rules for such an identification.

Yessenin-Volpin's work was influential; or perhaps one should say that it \emph{registered}. 
\cite{Parikh:Feasibility} refers to Yessenin-Volpin's claims.
Another key inspiration for Parikh's work seems to be \cite{Bernays:OnPlatonism}.\footnote{Parikh's work will likely be analyzed in other contributions to this volume.}
As already remarked, Yessenin-Volpin's work was also taken into account in the Czech school of mathematical logic, 
already in its inception with Rieger.
Recently, a major study in strict finitism has been carried out in \cite{Dean:strictFinitismFeasibilitySorites}.
Among its many contributions, the paper lays out in detail a suitable framework
(of weak two-sorted arithmetic) which makes it possible to formulate a statement about an existence of a bijection between numbers in unary notation
and numbers in binary notation, and exhibits an axiomatic theory which does not derive the  statement of the existence of a bijection.

\section{Dummett's analysis of strict finitism}

It is my intent to use the analysis of strict finitism
laid out in the paper \emph{Wang's paradox}, \cite{Dummett:WangParadox}, to mediate an interpretation of  Yessenin-Volpin's feasible numbers by means of countable semisets in the Alternative Set Theory.
Dummett's influential paper was published in a special issue of \emph{Synthese} dedicated to vagueness.
He does not purport to study feasibility; 
the vital problem in his work is to assess a lingering concern, going back to \cite{Bernays:OnPlatonism}, that 
 the arguments supporting intuitionism, especially if taken literally, are similar to those supporting strict finitism,
 which leaves an uncomfortably narrow space for intuitionism between platonism on its one side, and strict finitism
 on its other. This concern is eventually dismissed in the paper.
 
No citations are provided by Dummett, but \cite{Dean:strictFinitismFeasibilitySorites} 
refers to  \cite{Yessenin-Volpin:ProgrammeUltra} 
and \cite{Yessenin-Volpin:UltraIntCriticism1970} as Dummett's sources.\footnote{The endnote of \cite{Dummett:WangParadox} specifies that it records a lecture delivered in 1970;
and the two Yessenin-Volpin's works that \cite{Dummett:WangParadox} has been reconstructed to work from were presented at conferences in 1959 and in 1968 respectively.}

Dummett was working from a position of familiarity with Wittgenstein's work \emph{Remarks on the Foundations of Mathematics}, and he published his
appraisal of it in \cite{Dummett:WittgPhilMath}. It is also well known that Wittgenstein's \emph{Remarks}
were influenced by Hilbert's work---see especially \cite{Mulholzer:ProofMustBeSurveyable}
and \cite{Floyd:SurveyabilityHilbertWittgTuring}.
I am not aware of a similarly detailed analysis of   Yessenin-Volpin's work on this score.
 \cite{Dummett:WittgPhilMath} makes a cursory reference to \cite{Bernays:OnPlatonism} and writes:
\textquote[pp.~343--344]{As presented by Bernays, it [i.e., strict finitism] would consist in concentrating on
practical rather than on theoretical possibility. I have tried to
explain how for Wittgenstein this is not the correct way in which
to draw the contrast.} 
In \emph{Wang's paradox}, on the other hand, it is the emphasis on practice that is the distinctive trait of 
strict finitism as attributed to Yessenin-Volpin by Dummett.

\cite{Dummett:WangParadox} sets off by admitting that the strict finitist position 
admits \emph{multiple} totalities of natural numbers.
\begin{displayquote}[p.~302]
For by `natural number' must be understood a number
which we are in practice capable of representing. Clearly, capacity to
represent a natural number is relative to the notation allowed, and so the
single infinite totality of natural numbers, actual on the platonist view,
potential on the traditional constructivist view, but equally unique and
determinate on both, gives way to a multiplicity of totalities, each defined
by a particular notation for the natural numbers.
\end{displayquote} 
But quite early into his essay, Dummett goes on to consider a single totality
of natural numbers, opting for a notation (exemplified by a positional notation in an arbitrary integer base greater than 1) 
with the property that if a number $n$ can be represented in this notation, so can all the numbers less than $n$ be.
This seems to dismiss exponential functions as a general tool, although it does not totally exclude
the value  $10^{12}$ provided by Yessenin-Volpin from consideration, as one can hope to represent it
using just unary notation and multiplication and still get the sense across.
But Yessenin-Volpin preferred to consider multiple totalities that may perhaps be mapped onto one another in a particular manner and under conditions.\footnote{Cf.~the footnote on p.~9 of \cite{Yessenin-Volpin:UltraIntCriticism1970}: \textquote{Usually this is expressed by 
 `if $m \simeq n$ then $m' \simeq n'$. But such expressions overlook
the fact that $m' \simeq n'$ can appear only after $m \simeq n$ has appeared. 
This negligence may lead to the platonistic representation of a realm of all 
equivalences of the kind $m \simeq n$. Of course I have to reject this, [\dots]}}

Two key notions are introduced by Dummett: a \emph{weakly infinite} totality 
and a \emph{weakly finite} one. Coherence of the strict finitist position is then made to depend on
the existence of a totality that is both weakly infinite and weakly finite, as an extension of a vague predicate:
\textquote[p.~313]{As Esenin-Volpin
in effect points out, such totalities---those characterised as the extensions
of vague predicates---can be both weakly finite and weakly infinite.}
It is in the light of this interpretation of Yessenin-Volpin's intent that Dummett finally writes strict finitism off:
\textquote [p.~323]{As far as strict finitism is
concerned, common sense is vindicated: there are no totalities which are
both weakly finite and weakly infinite, and strict finitism is therefore an
untenable position.}

In turning to Dummett for a mediation, my intent is to take the notions he brought into the discussion 
and try to reinterpret them again in the Alternative Set Theory. 
Various motivations can be found for taking this step:
for example, \cite{Dummett:WangParadox} enjoys a lasting influence and the notions introduced therein
are more likely to be known to readers than the earlier setup of Yessenin-Volpin's work, 
or indeed the setup of Vopěnka's  Alternative Set Theory. 
And my \emph{main} reason for the attempt at a mediation is also one of familiarity: 
 the relative ease with which the Alternative Set Theory can capture a talk of weakly infinite, 
 weakly finite totalities, provided that a \emph{totality} is understood as a \emph{collection}.
In other words, the relevance of Dummett's work consists in shifting the discussion to a more familiar ground.
I do not intend to suggest  any direct influence.
 \cite{Vopenka:Teubner} does not refer to Dummett's work. 
 The relevance
 consists in Dummett highlighting the extensional perspective and \emph{a fortiori}, the perspective of collections and well-orders of them, an immediate effect of which is to render the subject matter under discussion amenable to the technical means and methods of set theory, inclusive of the Alternative Set Theory.

According to Dummett, a totality is
\begin{itemize}
\item  \emph{weakly infinite} if there exists a 
well-order of it with no last member; and 
\item \emph{weakly finite} if, for some finite ordinal $n$, there exists a 
well-order of it with no $n$th member. 
\end{itemize}
Dummett's running example---taken from \cite{Yessenin-Volpin:ProgrammeUltra}, p.203---is the
\emph{totality of heartbeats in my childhood}, ordered by occurrence in time.
Dummett tells us that \textquote{such a totality is weakly infinite, according to Esenin-Volpin: 
for every heartbeat in my childhood, I was still in my childhood
when my next heartbeat occurred.} 
On the other hand one can present a number $n$ (read: a closed arithmetical term, such as
$10^{12}$) such that the totality does not contain an $n$th member. 
Dummett also mentions what it would mean for a totality to be \emph{strongly finite}: 
namely to have a last member,  \textquote{like the set of heartbeats in my whole life}  
and \emph{strongly infinite}: not finitely bounded, \textquote{like the set of natural numbers, as ordinarily conceived, or, possibly, the set of heartbeats of my descendants.}
A totality with one of the last two attributes is \emph{determinate};
and one may succeed in embedding a vague, weakly infinite and weakly finite totality into a larger, determinate one.\footnote{\textquote{Hence, if induction is
attempted in respect of a vague predicate which in fact determines a
proper initial segment, which is both weakly finite and weakly infinite,
of a larger determinate totality, the premisses of the induction will both
be true but the conclusion will be false.} (p.~313)}

Let us pause at these examples for a heartbeat or two, to risk a pun.
Dummett's definition calls for a totality, which  I  read as a collection.
In the key part of his paper, where Dummett makes a commitment
to a \emph{bona fide} analysis of Yessenin-Volpin's viewpoint, we read that 
\begin{displayquote}[pp.~312--313]
Then we should normally say that a weakly finite totality could not also be weakly infinite. If we hold to this
view, we cannot take vagueness seriously. A vague expression will, in
other words, be one of which we have only partially specified a sense;
and to a vague predicate there will therefore not correspond any specific
totality as its extension, but just as many as would be the extensions of
all the acceptable sharpenings of the predicate. 
But to take vagueness
seriously is to suppose that a vague expression may have a completely
specific, albeit vague, sense; and therefore there will be a single specific
totality which is the extension of a vague predicate.
\end{displayquote}
The demonstrated worry (beyond the scope of \emph{Wang's Paradox} in the context of Dummett's work), is one about extensions of vague predicates---by his own admission, a legacy from Frege. 

On top of considering a totality, Dummett suggests to consider one or more well orders of it.
So  if one is looking for a totality that is weakly infinite and weakly finite,  the type of entities
searched for are (a) a collection, say $C$; (b) a well order of it, say $WO_1$, with no last element; 
and (c) a(nother) well-order $WO_2$ of the collection $C$, plus a background ordinal scale and technical tools 
that make it possible to express that $WO_2$ has no $n$-th member on this scale. 
The ambient theory in which these considerations take place thus allows a reference to finitude: 
namely a \textquote{finite} ordinal $n$. 
A few difficulties tied to such a requirement---of sustaining a suitable notion 
of finitude in the ambient theory, governed by logical principles aligning with the foundational position in question--- were demonstrated in Section 5 of \cite{Wright:StrictFinitism}. 
When Dummett speaks about \emph{membership} in an ordered totality, he refers not just to being an element of the collection, but rather to a \emph{position} in the well order of the collection.

Dummett ultimately rejects as incoherent the views advanced by Yessenin-Volpin, when the latter are interpreted
as the suggestion that weakly infinite, weakly finite totalities exist. 
His analysis rests on the incoherence of observational predicates, so in particular, these cannot be governed
by a consistent logic, as had been hoped.
In turn, the adequacy of Dummett's analysis to a genuinely strict finitist view 
was called into question by many authors, although none of them seek to reinstate Yessenin-Volpin's views.
In fact, if \cite{Wright:StrictFinitism} admits anyone as a promotor of strict finitism, it is Wittgenstein, rather than Yessenin-Volpin 
(but see also \cite{Marion:PhD} for an in-depth study of Wittgenstein's \emph{Remarks}).
In \cite{Dean:strictFinitismFeasibilitySorites}, the author takes into account an explicit rejection of assumptions (of ``traditional mathematics'') on  Yessenin-Volpin's part that some of Dummett's analysis is based on.

\section{On Vopěnka's Alternative Set Theory}

The Alternative Set Theory\footnote{I shall capitalize the name `Alternative Set Theory' and frequently use the acronym `AST'. 
This may run contrary to the customs of the original proponents. 
As to capitalization, many important resources on the theory are only available in Czech or Slovak, where the standards of capitalization differ from those of English. Not least, there is the practice of limiting the capitalization to the first word a phrase that functions as a proper name. In some Czech resources, the collocation ``alternativní teorie množin'' (``alternative set theory'') is not capitalized at all. Resources in English, especially the canonical  \cite{Vopenka:Teubner}, 
also observe this custom, i.e., no capitalization. This was significant in that it earmarked the \emph{alternative set theory} as a broader movement, rather than just an axiomatic theory. The use of acronym AST follows in the footsteps of, e.g., \cite{Sochor:AST76}; the original tradition may have reserved it for a particular axiomatic system. For the purpose of this paper I will ignore the distinction.}
 was developed by Vopěnka and a group of collaborators during the 1970s and 1980s. Two different monographs authored by Vopěnka were published in this period: \emph{Mathematics in the Alternative Set Theory}---\cite{Vopenka:Teubner}---in 1979 in English and \emph{Úvod do matematiky v alternatívnej teórii množín}---\cite{Vopenka:Alfa}---in 1989 in Slovak.
The group of collaborators notably included Antonín Sochor, who investigated the metamathematics in a series of papers \cite{Sochor:ASTmeta-I,Sochor:ASTmeta-II,Sochor:ASTmeta-III} and (co-)authored numerous papers on various mathematical topics within the theory, as well as several motivational works \cite{Sochor:AST76,Sochor:AST_CST} and the survey \cite{Sochor:BasesAST}.  With a few exceptions, papers about the theory came from authors within Czechoslovakia, and the publications on the AST gradually abated after Vopěnka shifted his interest to other topics after the year 1989. The inspiration for the theory came from several directions\footnote{This is studied in \cite{Hanikova:ATMvKanonu}.}, notably \emph{The Theory of Semisets}, \cite{Vopenka-Hajek:Semisets} with which it shares the key \emph{semiset} notion, and Robinson's Nonstandard Analysis as in \cite{Robinson:NSAbook}. 

From a logician's point of view, the Alternative Set Theory can be described in a nutshell 
as a first-order axiomatic nonstandard theory of sets and classes, developed in classical logic, where the set fragment is just the theory of hereditarily finite sets\footnote{In particular, the negation of the usual axiom of infinity is a theorem.}, and the class layer has comprehension for non-normal formulas,  an axiom of  prolongation that forces nonstandardness, and with a well-order axiom for classes.
The term `nonstandard' is used in the sense of Skolem's nonstandard models of Peano arithmetic ($\pa$). 
There is also a  sense in which said logician will take the theory to be quite ordinary: 
just  a particular axiomatic theory in classical logic.  

\cite{Vopenka:Teubner}  says  that \textquote[(p.~17)]{our concept of set is similar to Cantor's, but all our sets are finite from Cantor's point of view. We shall not admit the fiction of actually infinite sets.} 
The axioms for sets in the Alternative Set Theory prove all statements of the theory $\zffin$ 
(obtained by flipping the axiom of infinity for its negation and using the $\in$-foundation schema instead of regularity).
The converse is also true: the theory $\zffin$ proves all the axioms in the set fragment of the AST.
In particular these axioms can be taken as one way of introducing the theory $\zffin$  (cf.~Chapter I, Section 1 of \cite{Vopenka:Teubner}).
Sets are sometimes said to be \emph{formally finite} by the original proponents of the AST.
One can moreover show (in $\zf$) that the AST extends the ZFfin conservatively; cf.~\cite{Sochor:ASTmeta-III}.

A \emph{semiset} is a subclass of a set, thus a particular type of class (also, each set is a class, so classes are a universal sort). This definition of semiset coincides with the one given earlier in  \emph{Theory of Semisets} (TSS).
Both theories block the usual way of disproving the existence of proper semisets (from the axiom of replacement).
The existence of semisets, in either theory, is justified informally to someone who  already more or less endorses
the existence of sets and classes, is accustomed to fundamental ideas of set theory, and is now 
disposed to give their attention to the problem whether or not there are proper semisets.
One may also ask whether proper semisets are guaranteed to exist; in the AST this is the case, by prolongation.

Both the TSS and the AST exemplify viable axiomatic systems where the distinction between sets and classes does not follow the limitation of size criterion. This is pointed out in \cite{Fletcher:Infinity},\footnote{Fletcher considers only the $\atm$, not the $\tss$.} with references to \cite{Hallett:CantorianSetTheory} and \cite{Mayberry:Foundations}.\footnote{This is remarked early on by \cite{Sochor:AST76}:
``infinity is to be found already among the (formally) finite sets''.} 
Instead, Vopěnka stressed that the set vs.~class distinction arises from distinctness.
This is a keynote theme in  \cite{Vopenka:AstPhil1991}. 
 
A class in the AST is  \emph{finite} provided each of its subclasses is a set, otherwise it is infinite.
Thus each finite class is a set. 
The existence of proper semisets entails the existence of infinite sets, and these in turn entail 
the existence of infinite natural numbers (each set has a 1-1 set bijection to a unique natural number).
Infinite natural numbers in the AST were called ``inaccessible'' in \cite{Sochor:AST76}.\footnote{Moreover, 
a variant of prolongation was called ``the axiom of extension'' and the class $\FN$ was referred to as ``absolute numbers''.}
Following the notation of \cite{Vopenka:Teubner}, $\hat\approx$ denotes set bijections 
whereas $\approx$ denotes class bijections.
A useful theorem characterizes infinite sets in the AST, 
saying that a set $x$ is infinite if and only if $x\approx x\cup\{ y \}$  for each $y\not\in x$
(p.~37).
Infinite sets in the theory may thus be said to allow a rudimentary version of Hilbert's hotel, if  
 the maps involved are class maps.
A variant of the theorem says that a set is infinite if and only if it admits a class bijection to one its proper subsets.

A class $X$ in the AST is countable provided that it is infinite and there is a relation $R$ which is a total
order on $X$ and for each $x\in X$, the segment $\{y\in X \mid R(y,x)$ is finite. 
Such an order type is called $\omega$.\footnote{Below, we will encounter $\omega$ again as the smallest
infinite ordinal number.}
Thus no finite class is countable. 
The prolongation axiom says that any countable function can be extended to a set function.
As a simple consequence, each countable class is a proper semiset:
indeed any countable class $C$ can be obtained as a proper initial segment of a totally set-ordered 
set $a$, using the fact that $C$ cannot have a last element in the order
which is a witness to its countability (cf.~\cite[p.~42]{Vopenka:Teubner}) whereas total set orders always have a last element (\cite[p.~32]{Vopenka:Teubner}).
Countable classes can thus be taken as representing a new cardinality ``inserted'' between finite classes and classes such as  $V$ (the universal class) or $\N$ (natural numbers), which are not countable. 

The book \cite{Vopenka:Teubner} develops natural numbers in the AST. 
Vopěnka isolates facts that can be stated and proved in the theory $\zffin$ of hereditarily finite sets (not using classes,
except as shortcuts such as $V$ or indeed $\N$, the class of natural numbers). 
The class $\N$ equipped with the   successor operation, ordinal
addition and multiplication and the usual order ($x\leq y$ iff $x\in y$ or $x=y$) 
yields the structure  $\mathds{N}$ which gives an interpretation of  $\pa$ in $\zffin$.
For each set $x$ there is a unique $\alpha\in\N$ such that $x\,\hat\approx\, \alpha$, by induction (p.~59).  
On top of that (now using the class language), 
the AST has \emph{finite} natural numbers $\FN$. These  form a proper initial segment of $\N$
closed under successors (a proper cut). 
The cut is proper since an infinite set exists by prolongation, 
it has a set-bijection to a natural number, and this number cannot be finite, 
whereas $\FN$ is precisely the class of finite natural numbers.
The class $\FN$ is again  closed under addition and multiplication 
and so,  equipped with the restrictions of the operations and of the order, a structure $\mathds{FN}$ 
is obtained that gives an  interpretation of  $\pa$.
To define exponentiation on $\N$, let $\alpha^\beta$  be the (unique) natural number $\gamma$ such that
$\gamma \,\hat\approx\, \{ f : \mathrm{dom}(f) = \beta \mbox{ and } \mathrm{rng(f)} \subseteq \alpha\}$.  
Then $\FN$ is closed under exponentiation (cf.~\cite[p.~62]{Vopenka:Teubner}.\footnote{The closure follows from the closure of the class finite sets under other operations, notably the 
powerset and the cartesian product (cf.~p.~36--37) in \cite{Vopenka:Teubner}.
Then for $\alpha,\beta\in\FN$, by separation from $P(P(P(\alpha\cup \beta))))$ one gets $\alpha^\beta$, a finite set.
This maps to a unique $\gamma\in\N$, which must be finite (theorems p.~36 apply). 
Hence $\alpha^\beta\in \FN$, a theorem with a proof spanning a few paragraphs in the book.})
The class $\FN$ is infinite, and a typical example of a countable proper semiset.

\section{Feasible numbers and the Alternative Set Theory}

 \cite{Fletcher:Infinity}  points Vopěnka's Alternative Set Theory out as the most comprehensive attempt known to him  to model feasibility as standardness, in a suitable theory expressive enough to refer to standard and nonstandard natural numbers. 
\begin{displayquote}[p.~562]{It turns out that in AST the feasibility predicate can
be defined in terms of classes: a natural number n is feasible iff all subclasses of $\{0,1,\dots,n\}$
are sets.} 
\end{displayquote}
He also mentions a key aberration from the earlier approach of Yessenin-Volpin: 
in the AST, the value of any closed arithmetical term can be shown to be finite, i.e., an element of $\FN$.
The AST proves for example that the Bernays number $67^{257^{729}}$ is an element of $\FN$.
Fletcher suggests to execute the proof of feasibility of $10^{1000}$ in steps, 
repeating $10^{1000}$ times the step $F(n) \to F(n')$ 
(where $F$ is the feasibility predicate and ${}'$ denotes successor function).
Vopěnka also remarks on lengths of proofs in connection to the concept 
of a \emph{witnessed universe}: this postulates
the existence of a proper semiset in \textquote{a certain concrete set}.\footnote{In the quotation, I use the ``Bernays number''
from \cite{Bernays:OnPlatonism}, while the original has a different example.} 
\begin{displayquote}
    The study of witnessed universes is more difficult than that of limit universes. The theory of witnessed universes is in fact
inconsistent in the classical sense. If $c$ is an entirely concrete set (say, the set of all natural numbers less than
 $67^{257^{729}}$, then it can be obtained in finitely many steps from the empty set by successive addition of single elements; 
 thus $c$ is finite. On the other hand if $c$ has a proper subsemiset then $c$ is infinite in our sense. 
 But our proof of the fact that $c$ is finite has itself infinitely many steps (in our sense); thus it is not
convincing. 
\end{displayquote}
But the AST can fashion a reasonably short proof  that  $67^{257^{729}}$ is finite.
This proof is  an instantiation of the proof of the general fact stating that  $\FN$ is closed under exponentiation. 
Obtaining a proof that $67^{257^{729}}$ is finite amounts to instantiating the universally valid theorem  and its proof
with the numbers $67$, $257$, and $729$ somehow represented in PA, using numerals and arithmetical operations, a reasonably short proof.\footnote{There are models of the AST where $\FN$ has nonstandard elements. \cite{Sochor:BasesAST} comments as follows: \textquote{An interesting question is the relation between metamathematical natural numbers and elements of $\FN$. Evidently (formalization of) every metamathematical number is an element of $\FN$ however from the purely logical point of view $\FN$ can be `longer' [\dots] FN represents our smallest horizon and the system of metamathematical natural numbers forms a horizon and thus it is convenient to accept the (hardly formalizable) principle that $\FN$ formalizes exactly the system of metamathematical natural numbers.} (p.~66). The work \cite{Pudlak-Sochor:ModelsAST} presents a characterization of such 
cuts in countable models of $\pa$.} 
While witnessed universes are alluded to in \cite{Vopenka:Teubner}, to my knowledge 
no study of witnessed universes has been undertaken by Vopěnka or other researchers in the group. 

So in  broad strokes, to outline  whether and how the AST models feasible numbers, one 
might say they are modelled by the  proper cut $\FN$ on $\N$, 
definable   in the theory. 
The schema of existence of classes guarantees that $\FN$ exists;
the axiom of prolongation guarantees that $\FN$ is distinct from $\N$ (see above).

The term \textquote{feasible} is nowhere used in \cite{Vopenka:Teubner}.
On the other hand, \textquote{feasible}, as it occurs in the line of research
beginning with \cite{Yessenin-Volpin:UltraIntCriticism1970} (which is in English) and continuing with \cite{Parikh:Feasibility} was familiar to members of the Prague set theory seminar, given that it is used in \cite{Hajek:whySemisets} which reinterprets Parikh's work on almost consistent theories in the TSS. 
So taking also into account that Hájek was the translator of \cite{Vopenka:Teubner} into English,
the term might have found its way into the book.
One might speculate that, since the French cognate of \textquote{feasible} in  \cite{Yessenin-Volpin:ProgrammeUltra} is \textquote{réalisable}, at least a cognate term  indeed occurs in \cite{Vopenka:Teubner}, in the opening
passages, highly critical of classical set theory:
\begin{displayquote}[p.~9]
Cantor set theory is responsible for this detrimental growth
of mathematics; on the other hand, it imposed limits for mathematics that cannot be surpassed easily.
All structures studied by mathematics are apriori completed and rigid, 
and the mathematician's role is merely that of an observer describing them. 
This is why mathematicians are so helpless in grasping essentially
inexact things such as realizability, the relation of continuous
and discrete, and so on.
\end{displayquote}
This consideration could be pursued down a realm of translations of notions and putative shifts of their meaning, 
and while this might have some interest, I will leave it for another occasion.

There is a line of historical research carried out by Švejdar and in particular, outlined in \cite{Svejdar:VopenkaHajek}
and continued in \cite{Svejdar:VopenkaEarlySetTheory}.
This work  highlights the importance of the result of \cite{Vopenka-Hajek:ExistenceGenerelized}
where the authors show that G\"odel--Bernays set theory $\gb$ has proper definable cuts 
(whereas there are no proper definable cuts in $\zf$). 
The historical work traces the statement to  Vopěnka's habilitation thesis from 1964,
so he and Hájek probably knew this for some time prior to the publication of \cite{Vopenka-Hajek:ExistenceGenerelized}.
It is perhaps natural to wonder whether this insight may have played a role in the evolution of the AST.

This might seem to be the right moment to declare the case closed. 
It has been established that one distinctive feature of feasible numbers as rendered by Yessenin-Volpin, namely the unfeasibility of $10^{12}$, is not faithfully mirrored by the purported AST model.
On this score, Vopěnka might point out that the considerations above refer (only) to a \textquote{limit universe}, where there are indeed no proper subclasses  (i.e., subsemisets) of natural numbers given by closed arithmetical terms, because the limit universe
is an idealization. The examples in his book, on the other hand, pertain to witnessed universes.

But closing the case would leave a natural question unanswered, in connection to Dummett's weakly infinite, weakly finite totalities. Are there  proper semisets, as introduced in the AST, 
that can be shown to be both weakly infinite and weakly finite?
This question is naturally motivated by the fact that sets in the AST are 
\textquote{finite in the classical sense}, or in other words, \textquote{formally finite}. 
It is tempting to play the card of formal finitude of each natural number $\alpha\in\N$ and then,  
picking an $\alpha\in\N\setminus\FN$ to act in the role of the finite ordinal $n$ required by Dummett's definition of a weakly finite totality, to try to show that $\FN$ is both weakly infinite and weakly finite,  
provided only a well-order of some superclass of $\FN$ can be identified in which $\FN$ has no $n$-th member. 
After all, $\FN$ is \emph{bounded} in a suitable sense of the term.
Moreover, declarations to this effect were floated by the original proponents of the AST:  
\begin{displayquote}[\cite{Sochor:AST_CST}]
The classical set
theory puts infinity only `beyond' finite sets while the
alternative set theory puts it also `among' formally finite
sets, since there are sets containing infinite subclasses.
\end{displayquote}
Last but not least, one can view this route of accommodating feasible numbers in the AST
as parallel to Yessenin-Volpin's way of superposing the series of feasible numbers onto the series
that contains the element $10^{12}$.

It is also worth mentioning on this score that Vopěnka repeatedly made claims to vagueness, albeit using 
 different terms such as \textquote{inexactness} or \textquote{indefiniteness}:
\begin{displayquote}[\cite{Vopenka:AstPhil1991}]
Natural infinity is an abstraction from path leading to the boundaries of indefinite
and blurred phenomena. Since we understand various phenomena exactly through
grasping their boundaries (compare lat. `definitio'), the study of natural infinity is
also a possible foundation for the science of indefiniteness (p.~122).
\end{displayquote}
Some motivating examples for semisets are provided already in \cite{Vopenka:Teubner} that may be read as a plan to develop 
a mathematical theory of vagueness. One of his examples considers evolution of species, 
the successive generations that subsist between an ape on the one hand, and Darwin on the other. 
The sorites is referred to: \textquote[\cite{Vopenka:Teubner}]{Examples of proper semisets have been known for a long time but they were held for anomalies, as e.g. the `bald man paradox'.
But we meet proper semisets whenever in considering a property of
some objects we emphasize its intension rather than its extension. (p.~34)} 
 \cite{Dean:strictFinitismFeasibilitySorites} remarks that Vopěnka's proposal to attempt to capture vagueness by nonstandardness
 is exceptional in this respect.

\subsection{An attempt at weakly infinite, weakly finite totalities in the AST}

It is  fair to say in advance that the outcome of this subsection will be negative:  
 Dummett's weakly infinite, weakly finite totalities cannot be found in the AST.
The result as presented below concerns the AST and an attempt to reinterpret Dummett's notions in it, and this incurs granting the means of reasoning of classical logic (the underlying logic of the AST). 
It therefore does not say anything about the possibility of weakly infinite, weakly finite totalities in a more general setting, possibly not governed by classical logic. 
(Dummett's considerations, it seems fair to assume, were intended for a constructive reasoning.)

My aim here is to indicate the major obstacles.
As already stated, a first step of the intended reconstruction consists in reading a \emph{totality} as 
a \emph{class}.
So instead of a weakly infinite totality, 
I will henceforward take weakly infinite class, and likewise
in the case of a weak finitude.

Classes are often taken as the universal sort in the AST (cf.~\cite{Sochor:ASTmeta-I}); alternatively the theory can be taken as two-sorted with sets and classes but then one needs to say that each set is also a class (cf.~\cite[p.~28]{Vopenka:Teubner}). Either way, 
it is quite useful---vital even---to carefully consider what can be expressed using only the language of sets. 
Accordingly, I will talk about two \emph{views}:
\begin{itemize}
\item the \emph{set view} (statements that can be communicated by set formulas), and 
\item the \emph{class view} (statements that can be communicated by arbitrary formulas of the language of AST, i.e., using classes, and also using sets if desired).
\end{itemize}

Dummett's definitions of weakly infinite and weakly finite totalities (i.e., classes) 
each require a  well order to be supplied, along with the class under consideration.
Let us recall the definition of a well-order in the AST (\cite[p.~30]{Vopenka:Teubner}): 
$R$ is a well-order of $A$ if $R$ is a total order of $A$ and 
each non-empty \emph{subclass} of $A$ has a first element.           

To avoid a trivial solution, only  nonempty classes ought to be considered.
In the AST, one can freely read a  quantification into the definition.
A class $A$ can be declared to be weakly infinite provided there exists a well-order $R$ of it with no last member.
A well-order is binary relation on some domain, so it is a class in the AST. 
A class $B$ can be declared to be weakly finite provided there is a finite ordinal $n$ 
and a well order $S$ of $B$ with no $n$-th member.
For weak finitude, it is also useful to require that, 
if $S$ well-orders $B$ and there is an ordinal $n$ in the domain of $S$ such that
$B$ has no $n$-th member in $S$, then $n$ is an upper bound in the sense of $S$ to all elements of $B$.
Otherwise we might take any well order $S$ with a sufficiently large domain, 
 pick the least element of $S$, say $0$, and inject $B$ to $\mathrm{Dom}(S)\setminus 0$.

There is talk about ordinals,  in particular there is no $n$-th member of $B$ in the order $S$, 
so both $B$ and the domain of $S$ are required to be subsumed by ordinals.
I will simply assume that this is the case, for the class $\Omega$ of ordinals in the AST.\footnote{See \cite{Vopenka:Teubner}, Chapter II, Section 3 for the definition.} Namely weak finitude of $B$ is conditioned  as follows:
\begin{itemize}
\item $B \subseteq \Omega$ (B is  a class of ordinals or there is a bijection onto one, so $B$ inherits the order), and
\item there is an ordinal $\alpha$ such that $b<\alpha$ for each $b\in B$ ($B$ is bounded by $\alpha$), and
\item $\alpha$ is finite.
\end{itemize}

If $A$ is a nonempty class in the AST, 
$A$ is well ordered by $R$, and $A$ has no $R$-last element under $R$, then 
$A$ is infinite (on AST terms).
To see this, consider that $A$ is well ordered by $R$ and $\FN$ is well-ordered by $\leq$. 
By \cite[p.~31]{Vopenka:Teubner}, one of the two well-ordered classes is isomorphic to an initial segment
of the other. Suppose towards a contradiction that $A$ ordered by $R$ is isomorphic to a proper initial segment of $\FN$ ordered by $\leq$; let $B$ denote this proper initial segment of $\FN$, ordered by $\leq$.
Notice that $B$ has no last element in $\FN$ under $\leq$.
Let $c \in \FN \setminus B$, so $B$ is also a proper initial segment of  $c$ under $\leq$, with no last element. Now 
$c$ is a set, and $\leq$ on $c$ is a set order, so if $B$ is a set, it has a last element 
under $\leq$, a contradiction. So $B$ is not a set, hence $B$ is a proper semiset in $c$, 
a contradiction  since $c$ is finite.
Thus we conclude that $\FN$ well ordered by $\leq$ is isomorphic to an initial segment of  $A$ well ordered by $R$
and hence $A$ is infinite.

Can the class $A$, already known to be infinite,  
be also weakly finite? The quest for weakly finite AST classes has been restricted to 
the domain of ordinal numbers in the AST, i.e., to domain of the elements of $\Omega$. 
So the question is whether there is a class $A$ of ordinal numbers that, being weakly infinite, 
is also weakly finite.
We know that such a class is infinite on AST terms. 

Any class of ordinals in $\Omega$ is either cofinal in $\Omega$, or it is countable or finite. 
To see this, consider that $\Omega$ itself, as a class of ordinals, has the order type $\Omega$ (cf.~\cite[p.~69]{Vopenka:Teubner}): that is, for each $\beta \in \Omega$, the segment $\{\gamma \in \Omega \mid \gamma \leq \beta\}$ is countable or finite 
(cf.~the definition of order type $\Omega$, \cite[p.~50]{Vopenka:Teubner}).

Now our class $A$ under consideration cannot be cofinal in $\Omega$, 
since it is supposed to be bounded by an ordinal from above.
So $A$ is countable, and hence a proper semiset. 
For further consideration, I will fix one such countable proper semiset, namely the class $\FN$. 
Is this quite without a loss of generality?
 $\FN$ is a proper initial segment of $\Omega$, whereas other countable proper semisets of ordinals in $\Omega$
might not be. But should it turn out that $\FN$ could \emph{not} be bounded from above by a finite ordinal, 
then no more could any of the other countable ordinals be. So we can restrict attention to  $\FN$ without a loss 
of generality.

The following facts have emerged as relevant to the argument.

\begin{enumerate}
\item Any weakly infinite class (on Dummett's terms) is infinite (on AST terms).
\item $\FN$ is an example of a weakly infinite class, given it is well ordered by $\leq$.
\item $\FN$ is a proper class, in particular a countable semiset.
\item $\FN$ is not visible on the set view. 
\item $\FN$ is bounded by any infinite natural number  $\alpha \in \N\setminus \FN$. 
\item  No $\alpha\in \N\setminus \FN$ is  well-ordered by $\leq$ on the class view, 
  e.g., the class $\alpha\setminus\FN$ has no least element.
\item Each $\alpha\in\N$ is finite on the set view. 
\item Ordinal numbers in the AST form a subclass of $\N$ ($\Omega\subseteq \N$). 
The finite ordinals also belong to $\FN$. 
\item The class of ordinal numbers $\Omega$ is well ordered by (a restriction of) $\leq$. 
\item There is an ordinal $\omega \in \Omega$  such that
 $n\in\omega$ for each $n\in\FN$. The ordinal $\omega$ is the first element of $\Omega\setminus \FN$ (cf.~\cite[p.~69]{Vopenka:Teubner}).
\item One can find a well order of $\FN$ with no $\alpha$-th element: take the ordinal $\omega+1$ and let $\alpha \coloneqq \omega$.
\item  $\omega$ has infinitely many predecessors in $\Omega$, on AST terms.
\item The class $\Omega$ is  not visible on the set view. 
\item The set $\omega$ is visible on the set view. 
\end{enumerate}

Let us look at the ordinal $\omega$ in more detail. Since $\omega\in\N$, it is \textquote{formally finite}, 
as the proponents of AST used to say.
This formal finitude can be spelled out using the technical means of the AST and indeed, of its set fragment. 
To give one example, one can say that $\omega$ is Dedekind finite on the set view. 
$\omega$ as a set also meets all the conditions usually constituting the definition of an ordinal number;
but of course, on the set view, one cannot formulate the statement that $\omega\in\Omega$.

This list of facts seems to dash any hope that $\FN$ could be weakly infinite and weakly finite too.
Namely $\FN$ is infinite (2.) and not visible on the set view (4.). It is bounded by the ordinal $\omega$ (10.)
but not by a closed arithmetical term, and more importantly, the class of ordinals less than $\omega$ 
is infinite on the class view (12.). 
That is, it is infinite on the same view which makes $\FN$ visible.
So in other words, an attempt to use an ordinal such as $\omega+1$ to show that there is a well order 
in which $\FN$ has no $n$-th element, choosing $n=\omega$, can only be executed on the class view, 
 revealing that $\omega$ has infinitely many predecessors
  (i.e., the class $\FN=\{\alpha\in\Omega \mid \alpha \leq \omega\}$ is infinite). 

Are we really committed to use the same view (i.e., the class view) to argue weak finitude of $\omega$?
On the one hand, we have adopted the given definition of a well order in the AST, and this definition requires
all nonempty \emph{subclasses} to have least elements, so it requires the class view. 
This is why Vopěnka issues a warning saying that 
 \textquote{$\leq$ is no well-ordering oť $\N$} (\cite[p.~58]{Vopenka:Teubner}). 
 But on the other hand, surely, strictly on the set view (i.e., on the view of the theory $\zffin$), 
 $\N$ is well-ordered---namely, each nonempty \emph{subset} of $\N$ has a least element.
 
More broadly, perhaps  it would help to decline the mediation provided by Dummett's work at this point,
since it now appears that his notions requiring well orders did not contribute to addressing the  problem.
So here is a new attempt.  Fix an  $\alpha\in\N\setminus\FN$.
  All elements of $\FN$ (the feasible numbers) are bounded from above by $\alpha$ (5.),
so $\alpha$ has an infinite part. 
$\alpha$ is a (formally) finite natural number (7.) 
This finitude notion is different from the other finitude notion on which $\FN$ has been found to be infinite.
This attempt is thus carried out in full view of a dissociation of 
the meaning of \textquote{finite} from the meaning of \textquote{infinite}.
They are not complementary; to render this  in terms of their extensions, 
if applied to sets, their respective extensions might have a nonempty intersection and 
the union would not necessarily make up the whole universe of sets (i.e., universal class).  
In yet another words, \textquote{infinite} cannot be introduced in terms of \textquote{finite}. 
This is the price to be paid in order to accommodate feasible numbers in the AST.

Had this or similar approach to modelling feasibility been favoured by the Vopěnka group, 
it could  well have been put forward by them. 
As far as I can see, it never was.
Rather to the contrary, \cite{Vopenka:Teubner}  states in the Introduction that
the author has little to say on how his theory relates to Yessenin-Volpin's proposal:
\textquote{Motivations of certain approaches in our theory remind one
of the ideas of A.~S.~Esenin-Volpin; but the author feels unable
to present a more detailed analysis of these relations.}
Fletcher is also cautious in his formulation: \textquote{Now, if we alter the terminology, substituting the word ‘feasible’ for
‘standard’, then we can re-construe this as a theory of all the natural numbers, containing the
feasible numbers as a subsystem.}
It seems Vopěnka and Sochor never  put a foot wrong in selling their theory in this respect. 
In particular,  Sochor is careful to insert \textquote{formally} 
to his statement \textquote[\cite{Sochor:AST_CST}]{the
alternative set theory puts it [i.e., infinity] also `among' formally finite
sets}.

The so-called \textquote{formal finitude} is a misnomer in this sense.
Once a new notion of infinity (Vopěnka's natural infinity) has been introduced, 
the role of formal finitude is relegated to a set of tools that 
can  be used with sets. 
The set fragment of the $\atm$ is a weak theory (namely $\zffin$).
The powerset operation is available on sets (but not on classes uniformly). 
A total set order of a set has a last element. And so on.
In particular, the set we referred to as $\omega$ yields another set $P(\omega)$ on the set view.

On the other hand, once the requirement is restored that  \textquote{finite} and \textquote{infinite} be complementary notions, the ordinal $\omega$ in the AST cannot be considered finite for very similar reasons for which
the ordinal $\omega$ in classical set theory cannot be considered finite. 
Indeed Vopěnka remarks \cite[p.~68]{Vopenka:Teubner} that the ordinals in his theory
correspond to countable ordinals in Cantor's set theory. So in particular, 
$\omega$ in the AST would correspond in this way (i.e., as a well-ordering type) to the usual first infinite ordinal in $\zfc$ 
(in the sense of infinity used therein).

\section{Concluding remarks}

The approach taken above to modelling feasible numbers in the AST went via Dummett's weakly infinite, weakly finite 
totalities. There is a sense in which this attempt pulls heavily against the current.
I have earlier pointed out that ultrafinitism is viewed, and (at least with Yessenin-Volpin) also views itself,
as a departure from platonism in the direction of constructivism, only a more consistent one (it ``goes further'', as
Yessenin-Volpin stated). Since the AST is a theory in classical logic, it could be written off at any moment as
quite unfit for purpose on that score, even for accommodating Dummett's ideas, let alone Yessenin-Volpin's.

And yet, Vopěnka  alluded to Yessenin-Volpin for inspiration.
Moreover, one of the core problems laid out in Yessenin-Volpin's work, the problem of identity 
of (elements of) various number series he considers, can be addressed in an interesting way in the realm of classical logic. \cite{Dean:strictFinitismFeasibilitySorites} analyzes the relation between unary and positional 
notation in weak systems of two sorted arithmetic. A similar idea has also been included in  \cite{Pudlak:LogicalFoundations}.
Dean's work includes a detailed inquiry into the line of attacking  paradoxes of the sorites 
(or apparent paradoxes, were the attacks to be successful) by pointing out the essentially different 
tools that are used to indicate the \textquote{upper bound} from those that are  
available for \textquote{approaching it from below in steps}.  

I have laid out an argument for not taking countable semisets in the AST as examples of weakly infinite, weakly finite totalities. The argument is based on the definition of well-orders in the AST, which is presented on the class view (each nonempty \emph{class} has a least element in the order).
I have also tried to renounce a more direct route which uses two different notions of finitude to address the problem.
(But in Dummett's inquiry, of course, weakly infinite and weakly finite are not meant to be complementary notions.)
So in effect, I have been aiming at setting out some demarcation points 
regarding what the AST can contribute to the idea of feasible numbers and their models, 
if any, and what on the other hand it cannot contribute. 

One question set in the review of Yessenin-Volpin's work \cite{Yessenin-Volpin:ProgrammeUltra} by  \cite{Kreisel-Ehrenfeucht:ReviewYesseninVolpin} is the following:
\begin{displayquote}
    Can we find problematic properties of the intuitive notion of
feasible numbers, i.e., axioms satisfied by this notion, which [\dots]
can in turn be proved to be consistent in other conventional set theories?
\end{displayquote}
It has been one of my aims to enhance the idea  that the AST
models several important aspects of feasible numbers; namely, the class $\FN$ is 
a proper cut on $\N$, definable in the theory. The AST is relatively consistent, 
and in many respects  suitable for developing mathematics.

As a side effect, this paper may have called attention to 
other features of the AST and Vopěnka's work more generally. 
The novel notion of infinity in the AST, the AST its conception of sets and classes 
which does not follow the limitation of size principle, and the idea of proper definable cuts in $\gb$ 
are major examples of such features.

\medskip
\noindent {\bf Acknowledgements.}  
The work on this paper was supported by the Czech Science Foundation project 25-16489S.

\end{document}